\documentclass[12pt]{article}

\RequirePackage{amsthm,amsmath,amsfonts,amssymb,enumerate}

\RequirePackage[authoryear]{natbib}
\theoremstyle{plain}
\newtheorem{Lemma}{Lemma}
\newtheorem{Theorem}{Theorem}
\newtheorem{Proposition}{Proposition}
\newtheorem{Assumption}{Assumption}

\DeclareMathOperator*{\argmin}{arg\,min}
\def\R{{\mathbb{R}}}
\def\P{{\mathbb{P}}}
\newcommand\independent{\protect\mathpalette{\protect\independenT}{\perp}}
\def\independenT#1#2{\mathrel{\rlap{$#1#2$}\mkern2mu{#1#2}}}
\newcommand{\proofend}{$\hfill\Box{~}$}
\newenvironment{Proof}{\noindent {\em{\bf Proof.}}}{\proofend\\}    

\begin{document}
\begin{titlepage}

\title{High-dimensional nonconvex LASSO-type $M$-estimators }


\author{Jad Beyhum$^{1}$, François Portier$^{2}$ 
\\   \small $^{1}$ORSTAT, KU Leuven, jad.beyhum@gmail.com \\
        \small $^{2}$CREST, ENSAI, francois.portier@gmail.com \\ }

\maketitle

\begin{abstract}

This paper proposes a theory for $\ell_1$-norm penalized high-dimensional $M$-estimators, with nonconvex risk and unrestricted domain.
Under high-level conditions, the estimators are shown to attain the rate of convergence $s_0\sqrt{\log(nd)/n}$, where $s_0$ is the number of nonzero coefficients of the parameter of interest. Sufficient conditions {for our main assumptions} are {then} developed and {finally} used   {in several examples including} robust linear regression, binary classification and nonlinear least squares.

\end{abstract}

\smallskip

\noindent {{\large MSC 2020 subject classifications:} Primary 62F12; 62J02.}   \\

\smallskip

\noindent {{\large Key Words:} Lasso; High-dimensional regime; Nonconvexity; Unconstrained estimation.}   \\

\end{titlepage}



\section{Introduction}
Consider the standard statistical problem of estimating a parameter $\theta_0 \in \Theta \subset \mathbb R^d $, $d\geq 1$, defined as minimizing the unknown  \textit{true risk} $R: \Theta \to \mathbb R_{\geq 0}$ which is estimated by the \textit{empirical risk} $\widehat R : \Theta \to \mathbb R_{\geq 0} $ that depends on a random sample of size $n\geq 1$. The set $\Theta$ is called the parameter space and assumed to be convex (for simplicity). Motivated by \textit{large scale learning} applications, focus is on the \textit{high-dimensional} case in which the number $d$ of parameters is large relative to the sample size $n$.
To account for this situation, the asymptotic regime that shall be considered throughout the paper is the one of \textit{high-dimensional statistics} given by
\begin{align*}
n\to \infty  \qquad \text{and} \qquad d : = d_n \to \infty,
\end{align*} 
in which case, standard approaches, that directly minimize the empirical risk, are known to be inconsistent. The quantities $\theta_0,\Theta,R$ and $\widehat R$ implicitly depend on $n$, but we avoid to index them by $n$ to simplify the exposition. Reference textbooks dedicated to high-dimensional statistics include \cite{buhlmann2011statistics,giraud2015introduction,hastie2015statistical}.

A leading approach, that will be followed in this paper, is to \textit{regularize} the empirical risk by the $\ell_1$-norm of the parameters vector. In that, the estimate of $\theta_0$ is given by
\begin{equation}\label{est}
\widehat{\theta}\in\argmin_{\theta\in \Theta}  \, \{  \widehat{R}(\theta) +\lambda_n|\theta|_1\},
\end{equation}
where $\lambda_n>0$ is a penalty level that shall be chosen with respect to $n$.

Such a penalization approach, also referred to as the \textit{lasso}, has been successful in many cases such as \textit{linear regression} \citep{tibshirani1996regression, bickel2009simultaneous}, \textit{logistic regression} \citep{meier2008group} and \textit{Cox regression} \citep{tibshirani1997lasso,bradic2011regularization, huang2013oracle, kong2014non}. In presence of a \textit{sparsity structure} for $\theta_0$, i.e., when the number of nonzero coordinates $s_0:=s_{0n}$ of $\theta_0$ is small, the previous papers show that the \textit{lasso} method is reliable even in the challenging regime $ s_0\sqrt{\log(nd) / n }\to 0$. More specifically, results from the literature claim that the error $| \widehat{\theta} - \theta_0|_1$ is of order $s_0 \sqrt{ \log(nd) / n } $. While the results obtained for the three previous flagship examples, namely the \textit{linear}, the \textit{logistic} and the \textit{Cox regression}, are strong evidence of \textit{lasso}'s success, they all are developed for specific risk functions that are convex.  Extending the results to - still convex - but more general risk functions is the subject of recent work
such as \citep{van2008high,negahban2012unified}. In both cases the \textit{true risk} is globally (everywhere) convex and locally strongly convex.
Their main differences arise because \cite{van2008high} assumes that the \textit{true risk} is strongly convex on some neighborhood of the true parameter, while \cite{negahban2012unified} assumes that the \textit{empirical risk} is strongly convex on a cone of approximately sparse vectors. The latter assumption is called restricted strong convexity. 

More recently, the \textit{lasso} has been shown to be powerful in several
cases with nonconvex risk functions, e.g. \cite{yang2016sparse} for nonlinear least squares, \cite{stadler2010l} for mixture regression models, 
 \cite{loh2017statistical} for robust linear regression estimators and \cite{genetay2021high} for clustering. Even though these studies are carried out for specific estimates, only \textit{local} convexity on a small $\ell_2$ -ball around the true parameter is needed. {Note that \cite{stadler2010l} actually provides an oracle inequality for maximum likelihood estimators but does not obtain any rate of convergence on the estimation error.}
  Several papers \cite{wang2014optimal, loh2015regularized,mei2018landscape} propose high-level theories for nonconvex regularized high-dimensional $M$-estimators. In contrast to  \cite{yang2016sparse,loh2017statistical} (on specific applications), they require the strong convexity of the empirical risk on some sparse directions \cite{wang2014optimal} or on a cone \cite{loh2015regularized}.  
The empirical gradient and Hessian's behavior is investigated in \cite{mei2018landscape} but, concerning the asymptotic regime $ s_0\sqrt{\log(nd) / n }\to 0$,
no high-level result on the convergence of the estimator is given. 

The present paper establishes rates of convergence on  $| \widehat{\theta} - \theta_0|_1$ in the challenging regime $ s_0\sqrt{\log(nd) / n }\to 0$ without restrictive convexity assumptions.
The contributions can be summarized as follows:
\begin{enumerate}[(i)]
\item (generality) The proposed results are valid under a fairly general setting in that the risk function is not convex but only locally strongly convex in an $\ell_2$-ball around the true parameter. 
\item (interpretability and applicability) The results bear resemblance with well-known (low-dimensional) $M$-estimation theory \citep{newey1994large, van2000asymptotic, geer2000empirical} and can therefore be easily interpreted. We develop sufficient conditions for our high-level assumptions in order to simplify the application of the results.
\item (unrestricted domain) The proposed results do not require any restrictions on the parameter space. 
\end{enumerate}

As a secondary contribution, we apply our results to several examples including robust regression, binary regression, and nonlinear least squares. In each examples, the high-level results are easy to apply and the parameter space is $\Theta=\R^d$ illustrating the previous claims.

The fact that we allow the parameter space to be unrestricted may be surprising since the domain is restricted in \cite{stadler2010l, wang2014optimal, loh2015regularized, loh2017statistical, mei2018landscape}. This novel property is obtained through to a two-step technical argument. First, thanks to the penalization, we show that, regardless of $\Theta$ and with probability going to $1$, $\widehat \theta$ belongs to an $\ell_1$-ball $B$ with center $\theta_0$ and radius of order $\lambda_n^{-1}+|\theta_0|_1$. Second, the consistency of $\widehat \theta$ is obtained under an identification assumption on $R$ and a uniform convergence condition of $\widehat R$ on $B$. Because the radius of $B$ grows to infinity sufficiently slowly, the uniform convergence can be obtained in the applications of interest.

Note that a related but different problem is the one of computing the estimator \eqref{est}. Gradient descents algorithms usually converge to local minima of the objective function. Hence, $\widehat{\theta}$ may not be computable in practice. In the present paper, we do not consider this issue. Remark however that several papers \cite{wang2014optimal, loh2015regularized,yang2016sparse, loh2017statistical, mei2018landscape} treat both the optimization and statistical problems together by investigating the behavior of local minima of the function $\theta\in \Theta \mapsto \widehat R(\theta)+\lambda_n|\theta|_1$.

Another related line of work studies lasso-type estimators in the low-dimensional context where $d$ is fixed. In this setting, lasso-type estimators can be used for variable selection. For instance, \cite{fu2000asymptotics} develops an asymptotic theory for linear models in this framework and \cite{wang2013penalized} derives oracle properties in the more general case of possibly nonconvex semiparametric $M$-estimators. \\

\textbf{Outline.} In Section \ref{sec.hl}, we present the high-level results. Then, sufficient conditions for our high-level assumptions are stated in Section \ref{subsec.suff}. Next, we apply the results to three examples in Section \ref{sec.app}. Section \ref{sec.ccl} concludes the main text by discussing further research directions. The proofs of the high-level results and their sufficient conditions are in the Appendix. The results regarding the applications are proved in the supplement. \\

\textbf{Notations.} The notations $|\cdot|_1$, $|\cdot|_2$ and $|\cdot|_\infty$ correspond to the $\ell_1$, $\ell_2$ and sup norms, respectively. For a twice differentiable function $F:\R^K\mapsto\R$, $\nabla F$ is its gradient and $\nabla^2F$ its Hessian.

\section{High-level results}\label{sec.hl}

Let $\widehat R :\Theta\to \mathbb R_{\geq 0}$ be a random function and $R : \Theta\to \mathbb R_{\geq 0}$ be a function.
Consider $\hat \theta $ (resp. $\theta_0$) defined as a minimizer of  $\widehat{R}(\theta) +\lambda_n|\theta|_1$ (resp. $R$) over $\Theta$. In this section, the aim is to provide conditions on $\widehat{R}  $, $R$ and $\Theta$ to ensure certain convergence properties of $\hat \theta$ toward $\theta_0$. 

\subsection{Reduction of the parameter space} 

Thanks to the penalty term, we can show that, with probability going to $1$, $\widehat\theta$ belongs to an $\ell_1$-ball $B$ defined as
\begin{equation}\label{defB}B=\left\{\theta\in \Theta:\ |\theta|_1\le \lambda_n^{-1} (R(\theta_0) +1)   +  |\theta_0|_1\right\}.\end{equation}
This is formally claimed in the following proposition. 

\begin{Proposition}\label{prop:set}
We have that $\theta_0\in B$ and if $\widehat{R}(\theta_0) \to R(\theta_0)$, in probability, then $\widehat \theta \in B$ with probability going to $1$.
\end{Proposition}

This result is important because it allows in the mathematical development to restrict the attention to a smaller set $B$ included in the parameter set $\Theta$.  
The set $B$ has finite diameter (although its diameter can grow with $n$) while $\Theta$ could have infinite width. Hence, assumptions on the behaviour of $\widehat{R}$ are less demanding when restricted to $B$. This fact will be of good help when dealing with the applications. 

\subsection{Consistency} 

To obtain consistency, we make the following assumptions.

\begin{Assumption}\label{ident} For all $\eta>0$, there exists $\epsilon>0$ such that, for all $n\ge 1$,
$$\inf_{\theta\in \Theta, \, |\theta - \theta_0 |_2\geq \eta }\{R(\theta)-R(\theta_0)\}\ge \epsilon.$$
\end{Assumption}
This is an identification assumption restricting the shape of the true risk function. When $n$ is fixed, this condition holds if the risk is continuous and $\theta_0$ is its unique maximizer (the standard identification assumption in the literature of low-dimensional $M$-estimators). The specificity of the high-dimensional context is that we require this condition to be satisfied uniformly in $n$. In view of Proposition \ref{prop:set}, this assumption could be weakened by replacing $\Theta$ by $B$ but this does not bring much simplification because the set $B$ is intended to grow to $\Theta$ whenever $n$ is getting large.
It is also possible to relax Assumption \ref{ident} by letting $\epsilon$ go to $0$ with $n$. This could however prevent consistency if $\epsilon$ were to go to $0$ too quickly. {This has not been further investigated since Assumption \ref{ident} is valid in the applications considered in Section \ref{sec.app}.}\\


The second assumption ensures that the empirical risk converges uniformly to the true risk on $B$.
\begin{Assumption}\label{convLoss}
$\sup_{\theta\in B}\left| \widehat{R}(\theta)- R(\theta)\right|=o_P\left(1\right).$
\end{Assumption}
\noindent In the low-dimensional context, a similar condition is usually required on a compact set which does not depend on $n$. The main difference in the present context is that the radius of $B$ grows with $n$. 
\\

The following theorem states that the estimator is consistent in $\ell_2$-norm.
\begin{Theorem}\label{Consistency}
Under Assumptions \ref{ident} and \ref{convLoss}, if $\lambda_n|\theta_0|_1\to 0$, we have $| \widehat{\theta}-\theta_0 |_2=o_P(1)$.
\end{Theorem}
The result relies on the additional condition $\lambda_n|\theta_0|_1\to 0$, which, roughly speaking, means the added penalty term has only a negligible effect on the objective function evaluated at $\theta_0$.

\subsection{Rate of convergence} 

The following conditions are required to obtain a bound on the convergence rate of $\widehat \theta $ toward $\theta_0$.

\begin{Assumption}\label{Hessian}
There exist constants $\rho_*,\eta_*>0$ such that for all $n\ge 1$ and $\theta\in B, |\theta - \theta_{0} | _2 \le \eta_*$
$$R({\theta})-R(\theta_0)\ge\frac{\rho_*}{2}|{\theta}-\theta_0|_2^2.$$
\end{Assumption}
This is a local strong convexity assumption also imposed in the literature on nonconvex low-dimensional $M$-estimators. We stress that this condition is only imposed on an $\ell_2$-ball with radius fixed with $n$ (although the $\ell_2$-ball itself can change with $n$ since $\Theta$ and $\theta_0$ depends on $n$).  This condition does not require global convexity. Let $\mathcal{V}=\{\theta \in B,\ |\theta - \theta_{0} | _2 \le \eta_* \}$. A sufficient condition to obtain the previous assumption is to ask that $R$ is twice differentiable, $\theta_0$ is an interior point of $\Theta$, and the following eigenvalue property that for all $n\ge 1$,
$$\inf_{\theta \in \mathcal{V}} \rho_{\min}(\nabla^2R(\theta))\ge \rho_*,$$ where $\rho_{\min}(\cdot)$ is the minimal eigenvalue.
Indeed, as $\nabla R(\theta_0)=0$, by the second-order mean-value theorem, for all $n\geq 1$ and $\theta\in\mathcal{V}$, there exists $\tilde \theta\in\mathcal{V}$ such that
$$R(\theta)-R(\theta_0)=(\theta-\theta_0)^\top \frac{\nabla^2R(\tilde \theta)}{2}(\theta-\theta_0).$$\\

The last of our high-level conditions considers the difference between the empirical and the true risk 
\begin{align*}
\widehat \Delta(\theta )   = \widehat{R}({\theta})- R( {\theta}),
\end{align*}
and requires a certain convergence rate, $r_n$, for its increments.
\begin{Assumption}\label{convGradient}
There exist positive sequences $(r_n)_{n\geq 1} $ and $(\delta_n)_{n\geq 1}$ such that
\begin{align*}
 \lim_{n\to \infty } \P\left(\sup_{\theta\in\mathcal{V}}\frac{\left|\widehat{\Delta }( \theta)  -  \widehat{\Delta }( \theta_0)   \right|}{|\theta-\theta_0|_1\vee \delta_n}   \leq   r_n\right) = 1
\end{align*}
\end{Assumption}
A similar condition is also imposed in \cite{stadler2010l} for maximum likelihood estimators. In applications, $r_n$ and $\delta_n$ are typically of order $\sqrt{\log(nd)/n}$ and $\sqrt{\log(d)/n}$, respectively.

When the risk function is differentiable, the previous condition holds true as soon as, the gradient satisfies 
$\sup_{\theta\in \mathcal{V}}|\nabla \widehat{\Delta }(\theta) |_\infty\leq   r_n.$
Indeed in virtue of the mean value theorem, there exists $\bar \theta \in \mathcal{V}$ such that
\begin{align*} 
 \left|\widehat{\Delta }( \theta)  -  \widehat{\Delta }( \theta_0)    \right|= \left|\nabla \widehat{\Delta }(\bar \theta) ^\top ({\theta}-\theta_0)\right|  
  &\le  r_n \left|{\theta}-\theta_0\right|_1.
\end{align*}
As a result, Assumption \ref{convGradient} cares about the closeness (expressed trough $r_n$)  between the derivatives of the \textit{empirical risk} and the ones of the true risk. 

Remark also that the condition in  Assumption \ref{convGradient} depends on $\mathcal{V}$ which is itself defined through Assumption \ref{Hessian}. However, since $\mathcal{V}\subset B$, a stronger version of Assumption \ref{convGradient}  simply assumes \begin{equation}\label{as4'} \P\left(\sup_{\theta\in B}\frac{\left|\widehat{\Delta }( \theta)  -  \widehat{\Delta }( \theta_0)   \right|}{|\theta-\theta_0|_1\vee \delta_n}   \leq   r_n\right)\to 1,\end{equation}
where we stress that the supremum is taken on $B$ rather than on $\mathcal{V}$.  In Proposition \ref{sufficient} (see Section \ref{subsec.suff}), we provide sufficient conditions for the stronger result \eqref{as4'}.  The proof of Proposition \ref{sufficient} leverages empirical process theory. It avoids using the differentiability of the risk as outlined before. \\

Recall that $s_0$ is the number of non zero coordinates of $\theta_0$. We have the following Theorem.
\begin{Theorem}\label{Rate}
Under Assumption \ref{ident}, \ref{convLoss}, \ref{Hessian}, \ref{convGradient}, if $\lambda_n|\theta_0|_1\to 0$ and $\lambda_n \ge 2 r_n,$
with probability going to $1$, we have $$\left|\widehat{\theta}-\theta_0\right|_1\le \left(\frac{24}{\rho_*} s_0  r_n\right)\vee \delta_n.$$
\end{Theorem}
Since in the applications, $r_n$ and $\delta_n$ are of order $\sqrt{\log(nd)/n}$ and $\sqrt{\log(d)/n}$, respectively, Theorem \ref{Rate} gives us a rate of convergence of order $s_0\sqrt{\log(nd)/n}$, which is standard in high-dimensional statistics.

\section{Sufficient conditions} \label{subsec.suff}
In this section, we develop sufficient conditions for our high-level assumptions. They are leveraged to illustrate our theory with applications in Section \ref{sec.app}.

\subsection{Conditions on the true risk}

Two conditions are dealing with the function $R$, namely Conditions \ref{ident} and \ref{Hessian}. We here provide sufficient conditions, (i) and (ii) below, on the gradient $\nabla R$,  under which Conditions \ref{ident} and \ref{Hessian} are valid. They are based on the following proposition.

\begin{Proposition}\label{lemma:sufficient_cond}
Let $R: \Theta\to \mathbb R$ be differentiable and such that
\begin{itemize}
\item[(i)] For all $\theta\in  \Theta$, we have $\nabla R(\theta)^\top(\theta-\theta_0)\ge 0$.
\item[(ii)] For all $\gamma>0$, there exists $c(\gamma)>0$, decreasing in $\gamma$, such that, for all $n\ge 1$,
$$\inf_{\theta\in \Theta:\ |\theta-\theta_0|_2\le \gamma}\frac{\nabla R(\theta)^\top(\theta-\theta_0)}{|\theta-\theta_0|_2^2}\ge c(\gamma).$$
\end{itemize}
then for all $\theta\in \Theta$ and $ \eta >0$ such that $|\theta-\theta_{0}|_2 \ge \eta$,
\begin{equation}\label{into}R(\theta)-R(\theta_0)  \ge c(\eta) \frac{\eta^2}{2}. \end{equation}
This implies also that Assumptions \ref{ident} and \ref{Hessian} hold.
\end{Proposition}
Thanks to this proposition, only working on the function
  $\theta \mapsto \nabla R(\theta)^\top (\theta-\theta_0)$ is enough to obtain Assumptions \ref{ident} and \ref{Hessian}.

\subsection{Conditions on the empirical risk}\label{subsec:emp_risk}

Consider the standard regression setup where the goal is to predict $Y$, the response variable, with support $\mathcal{Y}\subset \R$, based on a random vector $X$ with support $\mathcal{X}\subset\R^d$. Let us use the notation
$$\mathcal{T}=\{x^\top\theta:\ x\in\mathcal{X}, \theta\in \Theta\}. $$ 
Interest is devoted to \textit{single index} types of risk defined as
\begin{equation}\label{Lg1}R(\theta)=E\left[\ell(X^\top \theta, Y)\right],\end{equation}
for all $\theta \in \Theta$,  where $\ell:\ (s,y)\in\mathcal{T}\times \mathcal{Y}\mapsto \R$.

Let $\{(Y_i,X_i)\}_{i=1}^n$ be an independent and identically distributed (i.i.d.) collection of random variables distributed as $(Y,X)$.  The estimate of $R$ is defined as 
\begin{equation}\label{Lg2}\widehat{R}(\theta)=\frac1n \sum_{i=1}^n \ell(X_i^\top\theta, Y_i).\end{equation}
We make the following assumption.
\begin{Assumption}\ 
\label{iid_single_index}  
\begin{itemize}
\item[(i)] There exists a constant $M_X>0$ such that $|x|_\infty \le M_X$ for all $x\in\mathcal{X}$.
\item[(ii)] There exists a constant $M_\ell>0$ such that 
$$\sup_{t\in\mathcal{T},y\in\mathcal{Y}}|\ell(t,y)| \le M_\ell.$$
\item[(iii)] There exists a constant $L>0$ such that 
$$\sup_{t,t'\in \mathcal{T},y\in\mathcal{Y}} \frac{|\ell(t,y)-\ell(t',y)|}{|t-t'|}\le L.$$
\end{itemize}
\end{Assumption}
The first two conditions stipulate that the features and the loss are bounded. The third condition imposes that the loss $\ell$ is Lipschitz with respect to its first argument uniformly in its second argument.

Let $|B|_1$ be the $\ell_1$-diameter of the set $B$ defined in \eqref{defB}, that is $$|B|_1=\sup_{\theta\in B}|\theta|_1.$$
We have the following proposition.

\begin{Proposition}\label{sufficient} Let $R$ and $\widehat R$ be defined as in \eqref{Lg1} and \eqref{Lg2}. Under Assumption \ref{iid_single_index}, if $\log(d)|B|_1^2n^{-1}\to 0$, then Assumptions \ref{convLoss} holds and property \eqref{as4'} is satisfied with 
$$\delta_n =\sqrt{\frac{\log(2d)}{n}},\ r_n= 
 16 LM_X \sqrt{   \log(4nd)  /n} .$$
\end{Proposition}
The fact that \eqref{as4'} is satisfied directly implies that Assumption \ref{convGradient} holds since $\mathcal{V}\subset B$. By definition of $B$, the assumption that $\log(d)|B|_1^2n^{-1}\to 0$ is both a condition on the rate of convergence to $0$ of $\lambda_n$ (which should not be too fast) and on the size of $|\theta_0|_1$ (which shall not be too large). Since $|\theta_0|_1$ and $s_0$ (the number of nonzero components of $\theta_0$) are strongly related, the latter can be interpreted as an assumption on the parsimony level $s_0$ which, roughly speaking, shall not exceed $\sqrt n $.

\section{Applications}\label{sec.app}
In this section, we show how to use the general results given previously to derive consistency results in specific applications, namely robust regression, binary classification and nonlinear least squares. 

\subsection{General regression setup}
Let us now introduce a regression framework that is similar to the one considered in Section \ref{subsec:emp_risk} but with some additional assumptions. This framework will be adopted in the three examples that follows. The response $Y$ has support $\mathcal{Y}\subset \R$, the covariates vector $X$ has support $\mathcal{X}\subset\R^d$. Let $\{(Y_i,X_i)\}_{i=1}^n$ be an i.i.d. collection of random variables with the same distribution as $(Y,X)$. The parameter space $\Theta$ is equal to $\R^d$. We further assume that $\theta_0\ne 0$ for simplicity. This condition ensures that the number $s_0$ of nonzero components of $\theta_0$ is strictly positive, which allows to simplify the statements of the rates of convergence.

The following assumption is made on the covariates vector.

\begin{Assumption}
\label{regX}There exists a constant $M_X$ such that for all $n\ge 1$, $\left|x\right|_\infty\le M_X$ for all $x\in\mathcal{X}$. The random vector $X$ has mean zero and is $M_X^2$ sub-Gaussian, that is $E[X]=0$ and
$E[e^{X^\top v}]\le e^{\frac{M_X^2|v|_2^2}{2}}$, for all $v\in\R^d$. There exists also $\rho_X>0$ such that for all $n\ge 1$, $\rho_{\min}(E[XX^\top])\ge \rho_X$.
\end{Assumption}
The fact that $X$ has bounded support allows to bound $X^\top\theta$ when $\theta$ lies in an $\ell_1$-ball. The sub-Gaussianity assumption ensures that $X^\top\theta$ remains small when $\theta$ lies in an $\ell_2$-ball. These two facts prevent  $X^\top\theta$ to take too large values which can have undesirable consequences on the estimation of the true risk and its shape. Then, the condition that $\rho_{\min}(E[XX^\top])\ge \rho_X$ is a classic identification assumption. The condition $E[X]=0$ can be easily avoided (at the cost of additional derivations) but is imposed for simplicity. Note that similar assumptions on the regressors are also imposed in \cite{mei2018landscape}.

\subsection{Robust regression}
We consider the following model:
$$\epsilon = Y - X^\top\theta_0 \quad \text{is such that} \quad   \epsilon \independent X \text{ and }  E[\epsilon]=0$$
and study robust estimators of the form 
$$\widehat{\theta}\in\argmin_{\theta\in \R^d} \frac1n \sum_{i=1}^n\rho(Y_i-X_i^\top \theta)+\lambda_n|\theta|_1,$$
where $\rho:\R\mapsto \R_+$ is some loss function and $\lambda_n>0$ is the penalty term.
This type of estimator falls into the class studied in Section \ref{subsec:emp_risk} with $\ell(t,y)=\rho(y-t)$. The distribution of $\epsilon$ is assumed to be independent of $n$ (unlike that of $X$, since the dimension $d$ of $X$ is allowed to grow with the sample size).
We make the following assumption. 
\begin{Assumption}
\label{NLS1}\ 
\begin{itemize}
\item[(i)] The mapping $\rho$ is continuously differentiable, $\rho'$ is odd, $\rho'(t)\ge 0$ for all $t\ge 0$ and there exists $M_\rho>0$, for all $t\in\R$, 
$$|\rho(t)|\vee|\rho'(t)|\le M_\rho.$$
\item[(ii)] The error term $\epsilon$ has a symmetric distribution, moreover, defining $$g:\ t\in\R\mapsto E[\rho'(t+\epsilon)],$$ we have $g(t)>0$ for all $t>0$, $g$ is differentiable at $0$ and $g'(0)>0$.
\end{itemize}
\end{Assumption}
This type of conditions can also be found in the robust regression example developed in \cite{mei2018landscape}. Condition (ii) is satisfied by Tukey's bisquare loss, which is usual in robust regression and given by
$$\rho_{\text{Tukey}}(t)=\left\{\begin{array}{cc}1-\left(1-\left(\frac{t}{t_0}\right)^2\right)^3&\text{for } |t|\le t_0\\
1& \text{for } |t|> t_0.
\end{array}\right.$$ Given that $\rho$ is odd and $\epsilon$ has a symmetric distribution, the condition that $g(t)>0$ for all $t>0$ holds when $\epsilon$ has a density which is strictly positive and decreasing on $\R_+$. 
We set
  $$\lambda_n=32M_\rho M_X\sqrt{\frac{\log(4nd)}{n}}.$$
We have the following Theorem.
\begin{Proposition}\label{TNLS1} Under Assumptions \ref{regX} and \ref{NLS1}, if $   \log(nd) |\theta_0|_1^2  n^{-1} \to 0$, we have 
$$\left|\widehat{\theta}-\theta_0\right|_1=O_P\left(s_0\sqrt{\frac{\log(nd )}{n}}\right).$$
\end{Proposition}
The additional condition $ |\theta_0|_1^2  \log(nd) / n \to 0$ is a sparsity condition, guaranteeing that $\lambda_n|\theta_1|\to 0$ and $\log(d)|B|_1^2n^{-1}\to 0$. In \cite{loh2017statistical} and \cite{mei2018landscape}, comparable results are obtained on robust regression estimators. On the one hand, in our case the parameter space is unrestricted while in \cite{loh2017statistical} (respectively, \cite{mei2018landscape})  it is limited to be an $\ell_1$-ball (respectively, $\ell_2$-ball). On the other hand, the rate of convergence derived by \cite{loh2017statistical} and \cite{mei2018landscape} applies to any local minimum of the estimation criterion while our rate only holds for the global minimum $\widehat\theta$. 

\subsection{Binary classification}\label{sec:binary_reg}

Suppose here that  $Y$ is binary, i.e., $\mathcal Y = \{0,1\}$ and consider the following model
\begin{align*}
\mathbb P ( Y=1 | X ) = \sigma(X^\top\theta_0 ),
\end{align*}
where $\sigma: \mathbb R \to  [0,1] $.
We define $R (\theta )  =  \mathbb E [ ( Y - \sigma ( X^\top\theta)  ) ^2 ]$ and let
\begin{align*}
\widehat \theta \in \argmin _{\theta\in \R^d }  \frac1n \sum_{i=1}^n  (Y_i - \sigma(X^\top\theta )  ) ^2 + \lambda _ n |\theta|_1.
\end{align*}
This estimator estimator belongs to the class studied in Section \ref{subsec:emp_risk} with $\ell(t,y)=(y-f(t))^2$. 
\begin{Assumption}
\label{assump:binary_reg}\ 
\begin{itemize}
\item[(i)]  The mapping $\sigma$ is differentiable, $\inf\limits_{|t| \leq s} \sigma'(t) > 0$ for all $s>0$ and  there exists a constant $M_\sigma>0$ such that $\sup\limits_{t\in \mathbb R}  \sigma'(t)  \le M_\sigma$.
\item[(ii)] There exists a constant $M_0>0$ such that, for all $n\geq 1$, $|\theta_0|_2\leq M_0 $. 
\end{itemize}
\end{Assumption}
Assumption \ref{assump:binary_reg} (i) means that $\sigma$ is strictly increasing and bounded. It imposes $\sigma'(t)\to 0$ as $t\to \pm\infty$. Such an assumption is, for instance, satisfied by the usual logistic function $\sigma(t)=(1+e^{-t})^{-1}$. We set 
$$\lambda_n=96M_\sigma M_X\sqrt{\frac{\log(4nd)}{n}}.$$ We have the following proposition.

\begin{Proposition}\label{prop:binary_reg}
Under Assumptions \ref{regX} and \ref{assump:binary_reg}, if $\log(nd)|\theta_0|_1^2 n^{-1}\to 0$,
we have 
\begin{align*}
\left|\widehat{\theta}-\theta_0\right|_1=O_P\left(s_0 \sqrt{\frac{\log(nd )}{n}}\right).
\end{align*}
\end{Proposition}

Remark that \cite{mei2018landscape} obtains a similar result under close assumptions. As for robust regression,  their rate of convergence holds for local minima of the objective function, but the parameter space is restricted to an $\ell_2$-ball. Note also that, thanks to Assumption \ref{assump:binary_reg} (ii) and the Cauchy-Schwarz inequality, one can show that $\log(nd)|\theta_0|_1 ^2 n^{-1}\to 0$ holds if $\log(nd)s_0n^{-1}\to 0$.

\subsection{Nonlinear least squares}\label{sec:nls} The last application studies the following model:
\begin{align*}
\epsilon = Y -  f(X^\top\theta_0) \quad \text{is such that} \quad  \epsilon| X \sim \mathcal{N}(0,\sigma^2),
\end{align*}
where $f: \mathbb R \to \R $ is a given function and $\sigma\ge 0$ is a constant.
Define $R(\theta )  =  \mathbb E [ ( Y - f ( X^\top\theta)  ) ^2 ]$ and let
\begin{align*}
\widehat \theta  \in \argmin _{\theta\in \R^d }  \frac1n \sum_{i=1}^n  (Y_i - f(X^\top\theta )  ) ^2 + \lambda _ n |\theta|_1.
\end{align*}

This case and that of binary classification are different because, here, the loss function may not be bounded nor Lipschitz since the support of $Y$ is no longer bounded. This prevents us from using Proposition \ref{sufficient}. To overcome this issue, we remark that 
\begin{equation}\label{decomp:nls} ( Y - \sigma ( X^\top\theta)  ) ^2 = (f(X^\top\theta)-f(X^\top\theta_0))^2+\epsilon^2-2\epsilon[f(X^\top\theta)-f(X^\top\theta_0)],\end{equation}
where the term $ (f(X^\top\theta)-f(X^\top\theta_0))^2$ is bounded and satisfies some Lipschitz property. The other terms can easily be handled using the Gaussian property on the distribution of $\epsilon$.

We make the following Assumption:
\begin{Assumption}
\label{assump:nls}\ 
\begin{itemize}
\item[(i)]  The mapping $f$ is differentiable, $\inf\limits_{|t| \leq s} f'(t) > 0$ for all $s>0$ and there exists a constant $M_f>0$, such that $\sup\limits_{t\in \mathbb R}  |f(t)|  <M_f$ and $\sup\limits_{t\in \mathbb R}  f'(t)  <M_f$.
\item[(ii)] There exists a constant $M_0>0$ such that, for all $n\geq 1$, $|\theta_0|_2\leq M_0 $.
\end{itemize}
\end{Assumption}
Assumption \ref{assump:nls} is similar to Assumption \ref{assump:binary_reg} with $\sigma$ replaced by $f$.
We have the following proposition.
\begin{Proposition}\label{prop:nls}
Under Assumptions \ref{regX} and \ref{assump:nls}, if $\log(nd)|\theta_0|_1^2n^{-1}\to 0$, there exists a constant $K>0$ such that, when $$\lambda_n\ge K\sqrt{\frac{\log(4nd)}{n}},$$
we have 
\begin{align*}
\left|\widehat{\theta}-\theta_0\right|_1=O_P\left(s_0\sqrt{\frac{\log(nd )}{n}}\right).
\end{align*}
\end{Proposition}
Note that also \cite{yang2016sparse} provides a rate a convergence for the local minimum of our estimation criterion in the context of high-dimensional nonlinear least squares. Their result is however obtained under a different set of assumptions. Indeed, they assume that there exists $a>0$ such that $f'(t)>a$ for all $t\in\R$. This is incompatible with our assumption that $f$ is bounded.

 \section{Conclusion}\label{sec.ccl}
Some additional research directions are of interest.
First, we could examine the variable selection properties of $\widehat \theta$. When a $\ell_1$-penalty is used, obtaining support recovery guarantees usually requires an incoherence assumption (see e.g. \cite{wainwright2009sharp} and references therein). 
Meanwhile, \cite{loh2017support} have shown that the incoherence condition can be avoided when nonconvex penalization schemes are used. Hence, attractive variable selections properties may be obtained in our general setup but with nonconvex regularizers (such as SCAD or MCP). Second, one may seek to study the prediction error in the present framework. In this case, the identification assumption may not be necessary and it should be possible to obtain oracle inequalities on the risk. \cite{stadler2010l} have obtained such results for maximum likelihood estimators. Finally, the behavior of semiparametric high-dimensional $M$-estimators could also be investigated.

\begin{appendix}
\section{Technical reminders } The results of this section are useful technical lemmas. They already appear in the supplementary material of \cite{beyhum2022extension} but are recalled to keep the paper self-contained.
\subsection{A Bound on the expectation} \begin{Lemma}\label{tech} Let $\{Z_i\}_{i=1}^n$ be i.i.d. mean zero $d$-dimensional random vectors such that $|Z_i|_\infty\le M$ almost surely for some constant $M>0$. Then, we have
$$E\left[\left|\frac1n\sum_{i=1}^n Z_i\right|_\infty\right]\le M\sqrt{\frac{2\log(2d)}{n}}.$$
\end{Lemma}
\begin{Proof} Take $v\in\R$ and $k\in\{1,\dots,d\}$. By Hoeffding's Lemma, we have $E[\exp(vZ_{ki})]\le \exp(v^2M^2/2)$. By independence of the $\{Z_i\}_{i=1}^n$, we obtain $E[\exp((v/n)\sum_{i=1}^nZ_{ki})]\le \exp(v^2M^2/(2n))$. For all $v>0$, this implies
\begin{align*}
E\left[\left|\frac1n \sum_{i=1}^nZ_i\right|_\infty\right]&= vE\left[\left|\frac1n \sum_{i=1}^nZ_i\right|_\infty/v\right]\\
&= v E\left[\log\left(\exp\left(\left|\frac1n \sum_{i=1}^nZ_i\right|_\infty/v\right)\right)\right]\\
&\le v \log\left(E\left[\exp\left(\left|\frac1n \sum_{i=1}^nZ_i\right|_\infty/v\right)\right]\right)\\
&\le v\log\left(\sum_{k=1}^d E\left[\exp\left(\frac{1}{n} \sum_{i=1}^nZ_{ki}/v\right)\right] + E\left[\exp\left(-\frac1n \sum_{i=1}^nZ_{ki}/v\right)\right]\right)\\
&\le v\log\left(2d\exp\left(\frac{M^2}{2nv^2}\right)\right)= v\left(\log(2d)+ \frac{M^2}{2nv^2}\right),
\end{align*}
where the first inequality is due to Jensen's inequality. Choosing $v=\sqrt{M^2/(2n\log(2d))}$ yields the result.
\end{Proof}

\subsection{Contraction theorem}
The following contraction theorem (Theorem 16.2 in \cite{van2016estimation}) will play an important role in our proofs. We now recall its statement for the sake of clarity.

\begin{Theorem}\label{contraction_th} Let $\{Z_i\}_{i=1}^n$ be a collection of random variables valued in $\mathcal Z$ and $\{\epsilon_i\}_{i=1}^n$ be a collection of independent Rademacher variables independent of $\{Z_i\}_{i=1}^n$. Let $\mathcal{F}$ be a class of functions defined on $\mathcal Z$ and valued in $\R$. Let $\rho : \mathbb R \times \mathcal Z \to \mathbb R $ be such that for all $(t,t') \in \mathbb R^2 $ and all $z\in \mathcal Z$,
$$|\rho( t, z) - \rho(t' , z)| \leq L_* |t- t'|, $$
for some $L_*>0$. Then, for all $f_*\in\mathcal{F}$, we have
$$E\left[\sup_{f\in\mathcal{F}}\left|\sum_{i=1}^n \epsilon_i(\rho(f(Z_i),Z_i )-\rho (f^*(Z_i), Z_i  ))\right|\right] \le2 L_*E\left[\sup_{f\in\mathcal{F}}\left|\sum_{i=1}^n \epsilon_i(f(Z_i)-f^*(Z_i))\right|\right].$$
\end{Theorem}
\begin{Proof}
Denote by $E_n$ the conditional expectation given $\{Z_i\}_{i=1}^n$. Applying Theorem 16.2 in \cite{van2016estimation}, we get
$$E_n\left[\sup_{f\in\mathcal{F}} \left|\sum_{i=1}^n \epsilon_i(\rho(f(Z_i),Z_i )-\rho(f^*(Z_i), Z_i  ))\right| \right] \le 2 L_*E_n\left[\sup_{f\in\mathcal{F}}\left|\sum_{i=1}^n \epsilon_i(f(Z_i)-f^*(Z_i))\right| \right].$$
We conclude using the law of iterated expectations.
\end{Proof}

\section{Proof of the results of Section \ref{sec.hl}} 
\subsection{Proof of Proposition \ref{prop:set}} 
The first statement is obvious because $R(\theta_0) $ is positive. Let us prove that $\widehat\theta\in B$ with probability approaching $1$. We have $\widehat{R}(\widehat{\theta})+\lambda_n|\widehat{\theta}|_1  \le\widehat{R}(\theta_0) +\lambda_n|\theta_0|_1$, which implies $\lambda_n|\widehat{\theta}|_1  \le\widehat{R}(\theta_0) +\lambda_n|\theta_0|_1$ because $\widehat{R}$ is positive. Therefore,  we obtain $\lambda_n|\widehat{\theta} |_1 =  R(\theta_0)   + o_P(1) + \lambda_n|\theta_0|_1 $. Hence with probability going to $1$, $\lambda_n|\widehat{\theta} |_1 \leq   ( R(\theta_0) + 1 ) +  \lambda_n|\theta_0|_1 $
 which yields the desired result.
\subsection{Proof of Theorem \ref{Consistency}} 
Let $\eta>0$. By Assumption \ref{ident}, there is $\epsilon>0 $ such that 
$\theta \in\Theta$ and $  |\theta-\theta_0|_2\ge \eta$ implies that $
R(\theta)-R(\theta_0)\ge \epsilon$. As a result, by the union bound,
\begin{align*}& P (|\widehat{\theta}-\theta_0|_2\ge \eta ) \leq P \left(  \{R (\widehat{\theta} )-R(\theta_0)\ge \epsilon\}\cap\{\widehat\theta\in B\} \right)+ P\left(\widehat\theta\notin B\right).
\end{align*}
By Proposition \ref{prop:set}, $\theta_0 \in B$ and so (by Assumption \ref{convLoss}), $ \widehat R(\theta_0 ) \to R(\theta_0)$, in probability. Then, invoking Proposition \ref{prop:set} again, the second term in the right-hand side goes to $0$. It remains to show that the first term  goes to $0$ as well. 

We have $R(\widehat{\theta})\geq R(\theta_0) $ and $\widehat{R}(\widehat{\theta})+\lambda_n|\widehat{\theta}|_1  \le\widehat{R}(\theta_0) +\lambda_n|\theta_0|_1$. It follows that, on the event $\{\widehat\theta\in B\}$,
\begin{align*}
0\leq R(\widehat{\theta})-R(\theta_0)&=[R(\widehat{\theta})-\widehat{R}(\widehat{\theta})] + [\widehat{R}(\widehat{\theta})-\widehat{R}(\theta_0)]+[\widehat{R}(\theta_0)-R(\theta_0)]\\
&\le 2\sup_{\theta\in B}\left|\widehat{R}(\theta)-R(\theta)\right| + \lambda_n(|\theta_0|_1-|\widehat{\theta}|_1 )\\
&\le 2\sup_{\theta\in B}\left|\widehat{R}(\theta)-R(\theta)\right| +\lambda_n |\theta_0|_1,
\end{align*}
where the second inequality is due to the fact that $\theta_0 \in B$. As a result
\begin{align*}
&P \left(  \{R (\widehat{\theta})-R(\theta_0)\ge \epsilon\}\cap\{\widehat\theta\in B\} \right)\\
&\leq P \left(  \{ 2\sup_{\theta\in B}\left|\widehat{R}(\theta)-R(\theta)\right| +\lambda_n |\theta_0|_1 \ge \epsilon\}\cap\{\widehat\theta\in B\} \right) \\
&\leq P \left(  2\sup_{\theta\in B}\left|\widehat{R}(\theta)-R(\theta)\right| +\lambda_n |\theta_0|_1 \ge \epsilon  \right) .
\end{align*}
In virtue of Assumption \ref{convLoss} and the fact that $\lambda_n|\theta_0|_1\to 0$,  the above term goes to $0$. 

\subsection{Proof of Theorem \ref{Rate}} 
For any set $\mathcal{E}\subset \Omega $, where $\Omega$ is the sample space of the probability space,  define $\mathcal{E}^c= \Omega\backslash \mathcal E $. Let $$\mathcal{A}= \left\{\left|\widehat{\theta}-\theta_0\right|_1> \left(\frac{24}{\rho_*} s_0  r_n\right)\vee \delta_n\right\},$$ $\mathcal{B}= \{\widehat\theta \in\mathcal{V}\}$ and $\mathcal{C} =\left\{\left| \widehat{\Delta }( \widehat \theta)  -  \widehat{\Delta }( \theta_0)   \right|\le r_n  | \widehat \theta- \theta_0| _ 1\right\}$. We have
\begin{align*}\notag P\left(\mathcal{A}\right) &= P\left(\mathcal{A}\cap \mathcal{B}\cap\mathcal{C}  \right)+P\left(\mathcal{A}\cap \left(\mathcal{B}\cap\mathcal{C}\right)^c  \right) \\
&= P\left(\mathcal{A}\cap \mathcal{B}\cap\mathcal{C}  \right)+P\left(\mathcal{A}\cap \left(\mathcal{B}^c\cup\mathcal{C}^c\right)  \right)\\
&\le P\left(\mathcal{A}\cap \mathcal{B}\cap\mathcal{C}  \right)+P\left(\mathcal{A}\cap \mathcal{B}^c  \right)+P\left(\mathcal{A}\cap \mathcal{C}^c  \right)\\
\notag  & := P_1+P_2+P_3,
\end{align*}
where the second inequality is due to the union bound.
It holds that $P_2\to 0$ because of Theorem \ref{Consistency}. Concerning $P_3$, since on $\mathcal{A}\cap \mathcal{C}^c $, it holds that 
$| \widehat \theta- \theta_0| _ 1 \vee \delta_n = | \widehat \theta- \theta_0| _ 1 $, we find that $\mathcal{A}\cap \mathcal{C}^c $ implies that 
\begin{align*}
\left| \widehat{\Delta }( \widehat \theta)  -  \widehat{\Delta }( \theta_0)   \right|>  r_n   ( | \widehat \theta- \theta_0| _ 1 \vee \delta_n) ,
\end{align*}
which, by Assumption \ref{convGradient}, has probability going to $0$. Hence, it suffices to show that $P_1\to 0$. We will show the even stronger result $P_1=0$.

Let $J=\text{Supp}(\theta_0)$. For a vector $v\in\R^d$, we denote by $v_J$ the vector in $\R^d$ such that $(v_J)_k=v_k$ for all $k\in J$ and $(v_J)_k=0$ otherwise. We also define $v_{J^c}=v-v_J$. 
Throughout the rest of the proof, we work on the event $\mathcal{B}\cap\mathcal{C}$.
By definition of $\widehat{\theta}$, we have
$$\widehat{R}(\widehat{\theta})+\lambda_n|\widehat{\theta}|_1\le \widehat{R}(\theta_0)+\lambda_n|\theta_0|_1.$$
Next, remark that
\begin{align*}|\theta_0|_1-|\widehat{\theta}|_1&=|\theta_0|_1-|\widehat{\theta}-\theta_0+\theta_0|_1\\
&=|\theta_0|_1-|(\widehat{\theta}-\theta_0)_J +\theta_0|_1- |(\widehat{\theta}-\theta_0)_{J^c}|_1\\
&\le|(\widehat{\theta}-\theta_0)_J|_1- |(\widehat{\theta}-\theta_0)_{J^c}|_1,
\end{align*} 
where we have just used that $|a|_1 - |b |_1 \leq |a-b|_1 $. Therefore, it holds that 
\begin{align*}
R(\widehat{\theta})-R(\theta_0)&  =   \left\{   \widehat{\Delta }(\theta_0) - \widehat{\Delta }(\widehat \theta ) \right\} +   ( \widehat R(\widehat{\theta})-  \widehat  R(\theta_0) ) \\
&\le   r_n   | \widehat \theta- \theta_0| _ 1   +  \lambda_n\left(|(\widehat{\theta}-\theta_0)_J|_1- |(\widehat{\theta}-\theta_0)_{J^c}|_1\right).
\end{align*}
By Assumption \ref{Hessian}, we have
\begin{equation*}
R(\widehat{\theta})-R(\theta_0)\ge\frac{\rho_*}{2}|\widehat{\theta}-\theta_0|_2^2\ge \frac{\rho_*}{2s_0}|(\widehat{\theta}-\theta_0)_J|_1^2. 
\end{equation*}
which implies 
\begin{align*}
\frac{\rho_*}{2s_0}|(\widehat{\theta}-\theta_0)_J|_1^2&\le r_n\left|\widehat{\theta}-\theta_0\right|_1 + \lambda_n\left(|(\widehat{\theta}-\theta_0)_J|_1- |(\widehat{\theta}-\theta_0)_{J^c}|_1\right)\\
& = r_n  \left( 3 |(\widehat{\theta}-\theta_0)_J|_1-  |(\widehat{\theta}-\theta_0)_{J^c}|_1\right)
\end{align*}
which yields that $3|(\widehat{\theta}-\theta_0)_J|_1\ge |(\widehat{\theta}-\theta_0)_{J^c}|_1$ and $|(\widehat{\theta}-\theta_0)_J |_1\le 6 s_0 r_n/\rho_*$. Hence, on the event $\mathcal{B}\cap\mathcal{C}$, we have $\left\{\left|\widehat{\theta}-\theta_0\right|_1\le \frac{24}{\rho_*} s_0  r_n\right\}$, which implies $P_1=0$.

\section{Proof of the results of Section \ref{subsec.suff}} 

\subsection{Proof of Proposition \ref{lemma:sufficient_cond}}

Take $\eta>0$ and $\theta\in \Theta,\ |\theta-\theta_0|_2\ge \eta$. Let us consider the mapping $t\in[0,r]\mapsto R(\theta_t)$, where $r=|\theta-\theta_0|_2$ and $\theta_t=\theta_0+tr^{-1}(\theta-\theta_0)$. This mapping is differentiable, with derivative $r^{-1}\nabla R(\theta_t)^\top(\theta-\theta_0)$. As a result, we have
\begin{align*}\notag R(\theta)-R(\theta_0) &=\int_0^r r^{-1}\nabla R(\theta_t)^\top(\theta-\theta_0)dt\\
\notag & = \int_0^r t^{-1} \nabla R(\theta_t)^\top(\theta_t-\theta_0)dt\\
&\ge  \int_0^\eta t^{-1}   \nabla R(\theta_t)^\top(\theta_t-\theta_0)dt,
\end{align*}
where the inequality results from (i).  By (ii) we have, for all $t\in[0,\eta]$, $$\nabla R(\theta_t)^\top(\theta_t-\theta_0)\ge c(\eta)|\theta_t-\theta_0|_2^2=c(\eta)t^2r^{-2}|\theta-\theta_0|_2^2=c(\eta) t^2.$$ Taking the integral yields \eqref{into}.

Let $\eta>0$. By taking the infinum over $|\theta-\theta_{0}|_2 \ge \eta$ in \eqref{into}, we obtain Condition \ref{ident}. To obtain Condition \ref{Hessian}, remark that \eqref{into} applied with $\eta=|\theta- \theta_0|$, leads to $R(\theta)-R(\theta_0)\ge c(|\theta-\theta_0|_2) {|\theta-\theta_0|_2^2} / {2}$ for all $\theta\in S$. Pick $\eta_*>0$. Since $c(\cdot)$ is decreasing, for any $\theta\in \Theta,\ |\theta-\theta_0|_2\le \eta_*$, we get $$R(\theta)-R(\theta_0)\ge c(\eta_*)\frac{|\theta-\theta_0|_2^2}{2},$$ which proves that Assumption \ref{Hessian} is satisfied too.

\subsection{Proof of Proposition \ref{sufficient}}

\subsubsection{An auxiliary result}
Recall that, for all $\theta\in \Theta$, we have $\widehat{\Delta}(\theta)=\widehat R(\theta)-R(\theta)$. For a bounded subset $C$ of $\Theta$, define $\mu_{C}= \sup_{\theta\in C} |\theta-\theta_0|_1$.
\begin{Lemma}\label{auxiliary}
Under the assumptions of Proposition \ref{sufficient}, for all bounded subsets $C$ of $\Theta$, we have 
$$\P\left(D \ge  2K\sqrt{\frac{2\log(2d)}{n}}\mu_C+K\mu_Ct\right)\le \exp\left(-\frac{nt^2}{8}\right) ,$$
where $K=2LM_X$ and $D =\sup\limits_{\theta\in C}\left|\widehat{\Delta}(\theta)-\widehat\Delta(\theta_0)\right|.$
\end{Lemma}
\begin{Proof}
Let $\{\epsilon_i\}_{i=1}^n$ be i.i.d. Rademacher variables independent of $\{(X_i,Y_i)\}_{i=1}^n$. Applying the symmetrization theorem (Theorem 16.1 in \cite{van2016estimation}), we obtain 
\begin{align*}E[D ]&\le 2E\left[\sup_{\theta\in C}\left|\frac1n\sum_{i=1}^n \epsilon_i[\ell(X_i^\top\theta, Y_i)-\ell(X_i^\top\theta_0, Y_i)]\right|\right].
\end{align*}
 The function $\ell ( \cdot, y) $ is Lipschitz on $\mathcal{T}$  by Assumption \ref{iid_single_index} (iii). It can be extended on $\mathbb R$ taking instead
  $ \ell(p(s) ,y)$ where $p(s)$ is the projection on  $\mathcal T_n$. This new function coincides with $\ell ( \cdot, y) $ and has the same Lipschitz constant. 
 Hence, we can apply Theorem \ref{contraction_th} to obtain
$$E[D ]\le 4L E\left[\sup_{\theta\in C}\left|\frac1n\sum_{i=1}^n \epsilon_iX_i^\top(\theta-\theta_0)\right|\right].$$
By H\" older's inequality, Lemma \ref{tech} and the definition of $\mu_C$, we get
\begin{align}\notag E[D ]&\le 4LE\left[\left|\frac1n\sum_{i=1}^n \epsilon_iX_i\right|_\infty\right] \sup_{\theta\in C} |\theta-\theta_0|_1\\
\label{bound_exp}&\le 2K\sqrt{\frac{2\log(2d)}{n}}\mu_C.
\end{align}
 
\noindent\textit{Concentrating $D$ around $E[D]$.} 
Because 
$$\left\{\ell(X_i^\top\theta, Y_i)-\ell(X_i^\top\theta_0, Y_i)-\{R(\theta)-R(\theta_0)\}\right\}_{i=1}^n$$ are i.i.d., bounded by $K\mu_C$ and mean zero, we can apply Massart's inequality (Theorem 16.4 in \cite{van2016estimation}). It leads to 
$$ \P\left(D \ge  KE[D ]+K\mu_Ct\right)\le \exp\left(-\frac{nt^2}{8}\right),$$
for all $t>0$. Using \eqref{bound_exp}, we obtain the result.
\end{Proof}

\subsubsection{Proof that Assumption \ref{convLoss} holds}\label{subsecconvloss}

Define $A =\sup\limits_{\theta\in B}\left|\widehat{\Delta}(\theta)\right|.$ Remark that $|A|\le D+\left|\widehat{\Delta}(\theta_0) \right|$, where $D =\sup\limits_{\theta\in B}\left|\widehat{\Delta}(\theta)-\widehat\Delta(\theta_0)\right|.$
We now show that both terms go to $0$ in probability. First, notice that the random variables $\{\ell(X_i^\top\theta_0, Y_i)-R(\theta_0)\}_{i=1}^n$ are unidimensional i.i.d. random variables with mean zero and bounded almost surely by $2M_\ell$. Hence, by Hoeffding's inequality, we have $ \left|\widehat{\Delta}(\theta_0) \right|=o_P(1).$ Second, applying Lemma \ref{auxiliary}, we obtain
$$ \P\left(D \ge  2K\sqrt{\frac{2\log(2d)}{n}}|B|_1 +K|B|_1t\right)\le \exp\left(-\frac{nt^2}{8}\right),$$
which yields (choosing $t = M/\sqrt n$ with $M$ large) that $D=o_P(1)$ (since $\log(d)|B|_1^2/n\to 0$). As a result, $|A|=o_P(1)$.

\subsubsection{Proof that property \eqref{as4'} holds}\label{subsecconvgradient}
Define the collection of rings $C_k=\{\theta \in B:\ 2^{k}\le |\theta-\theta_0|_1\le 2^{k+1}  \}$,  $k\in\mathbb{Z}$. Define
\begin{align*}
&A = \sup_{|\theta - \theta_0|_1 \leq \delta_ n } \left|\widehat{\Delta}(\theta)-\widehat\Delta(\theta_0)\right|  \\
& A_k =\sup\limits_{\theta\in C_k}\left|\widehat{\Delta}(\theta)-\widehat\Delta(\theta_0)\right|.
\end{align*}
Define $ u_n=\left\lfloor \frac{\log\left(\delta_n\right)}{\log(2)}\right\rfloor$ and $  v_n=\left\lceil \frac{\log\left(B  \right)}{\log(2)}\right\rceil$ and further assume that $\delta_n \leq B$ (the case $\delta_n > B$ is simpler as only $A$ needs to be bounded).  Since
  $\{ |\theta - \theta_0|_1 \leq \delta_ n \} \cup \{\cup_{k=u_n} ^{v_n} C_k\} $ covers the set $B$ and because $|\theta-\theta_0|_1\ge 2^k $ on $C_k$, it holds that
\begin{align*}
\frac{\left|\widehat{\Delta}(\theta)-\widehat\Delta(\theta_0)\right|}{|\theta-\theta_0|_1\vee \delta_n} \leq (A\delta_n^{-1} )\vee  (\max_{u_n\leq k  \leq v_n } A_k { 2^{-k} })   .
\end{align*}
We now focus separately on each of the two terms appearing in the above upper bound.
For the term in the right, by Lemma \ref{auxiliary}, we have
$$\P\left(A_k \ge  2K\sqrt{\frac{2\log(2d)}{n}}2^{k+1}+ K2^{k+1}t\right)\le \exp\left(-\frac{nt^2}{8}\right),$$
and the union bound yields
\begin{equation}\label{boundk}\P\left(\max_{ u_n \leq k \leq v_n} 2^{-k} A_k \ge  4K\sqrt{\frac{2\log(2d)}{n}}+ 2Kt\right)\le (v_n-u_n+1)\exp\left(-\frac{nt^2}{8}\right)\end{equation}
For the term in the left, 
by Lemma \ref{auxiliary}, it holds that
\begin{equation}\label{boundtilde}\P\left({A} \delta_n ^{-1} \ge  2K\sqrt{\frac{2\log(2d)}{n}}+ K t\right)\le \exp\left(-\frac{nt^2}{8}\right).\end{equation}
The union bound and the two inequalities \eqref{boundk} and \eqref{boundtilde} together gives that for all $t>0$,
$$\P\left(\sup_{\theta\in B}\frac{\left|\widehat{\Delta}(\theta)-\widehat\Delta(\theta_0)\right|}{|\theta-\theta_0|_1\vee \delta_n} \ge 4K\sqrt{\frac{2\log(2d)}{n}}+ 2Kt\right)\\
 \le (v_n-u_n+2) \exp\left(-\frac{nt^2}{8}\right)$$
Choosing $t_n=\sqrt{8\log(2n)/n}$ and using that $ \lceil x\rceil  \leq x + 1$, we obtain that
\begin{align*}
&\P\left(\sup_{\theta\in B}\frac{\left|\widehat{\Delta}(\theta)-\widehat\Delta(\theta_0)\right|}{|\theta-\theta_0|_1\vee \delta_n} \ge 
2K \left(\sqrt{\frac{ 8\log(2d)}{n}}+ t_n \right) \right) \le \frac{v_n - u_n  +2}{2n} \le \frac{\log_2(B) +3}{2n} .
\end{align*}
By assumption, it holds that $B^2 = o(n)$ which implies that the previous probability goes to $0$. 
Hence, taking $R =  8K \sqrt{   \log(4nd)  /n} $ we get 
$$\P\left(\sup_{\theta\in B}\frac{\left|\widehat{\Delta}(\theta) -\widehat\Delta(\theta_0)\right|}{|\theta-\theta_0|_1\vee \delta_n} \ge 
R\right)\le \P\left(\sup_{\theta\in B}\frac{\left|\widehat{\Delta}(\theta) -\widehat\Delta(\theta_0)\right|}{|\theta-\theta_0|_1\vee \delta_n} \ge 2K \left(
\sqrt{\frac{ 8\log(2d)}{n}}+ t_n\right)\right)$$
where the inequality follows from the fact that $2K\left(\sqrt{\frac{ 8\log(2d)}{n}}+ t_n\right) \leq r_n$ which comes from $\sqrt {a} +\sqrt {b} \leq \sqrt 2 (\sqrt{ a + b} )$ for all $a, b >0$. Since the previous upper bound goes to $0$ we have just obtained \eqref{as4'}.

\end{appendix}

\bibliographystyle{apalike} 
\bibliography{Paper}       


\end{document}


\title{Supplement to ``High-dimensional nonconvex lasso-type $M$-estimators''}


\author{
{\large Jad B\textsc{eyhum}}
\footnote{ORSTAT, KU Leuven and CREST, ENSAI}\\\texttt{\small jad.beyhum@gmail.com}
\and 
{\large Fran\c cois P\textsc{ortier}} \footnote{CREST, ENSAI.}
\\\texttt{\small francois.portier@gmail.com}
}

\date{}

\maketitle

The supplement contains some technical reminders (Section \ref{sec.tech}) and the proofs of Propositions 4, 5 and 6 of the main text (respectively in Sections \ref{sec.prop4}, \ref{sec.prop5}, \ref{sec.prop6}).

\setcounter{section}{0}
\section{Technical reminders on (sub-)Gaussian random variables}\label{sec.tech}

\begin{Lemma}\label{lem:subG2} 
Let $Z$ and $Z'$ be two sub-Gaussian random variables with variance parameter $v$ and $v'$, respectively. Then for all $\delta>0$, it holds
\begin{align*}
&\mathbb E \left[ Z^2 1_{| Z'|>  2 \sqrt { v'  \log(   4\sqrt 2  / \delta ) } }   \right] \leq  v  \delta.
\end{align*}

\end{Lemma}

\begin{Proof}
 Use Cauchy-Schwarz inequality to obtain $\mathbb E [ Z^2 1_{| Z'| > t}   ]\leq \mathbb E [ Z^4 ]^{1/2} \mathbb P ( | Z'| > t ) ^{1/2}  $. Use  Theorem 2.1 in \cite{boucheron2013concentration} to get $\mathbb E [ Z^4 ]^{1/2}  \leq 4v$. By definition of sub-Gaussian variables, $\mathbb P ( | Z'| > t ) \leq 2 \exp( -t^2 / 2 v') $. It follows that
\begin{align*}
\mathbb E [ Z^2 1_{| Z'| > t}   ] \leq 4 v \sqrt{2  \exp( -t^2 / 2 v')}=4 v \sqrt 2  \exp( -t^2 / 4 v'),
\end{align*}
and using $t = 2 \sqrt { v'  \log(  4 \sqrt 2 /  \delta ) } $ we obtain the statement of the lemma.
\end{Proof}

\begin{Lemma}\label{tech2} Let $\{\epsilon_i\}_{i=1}^n$ be i.i.d. $\mathcal{N}(0,\sigma^2)$ random variables for a some $\sigma\ge 0$ and $\{X_i\}_{i=1}^n$ be i.i.d. $d$-dimensional random vectors independent of $\{\epsilon_i\}_{i=1}^n$ . Assume that there exists a constant $M_X>0$ such that such that $|X_i|_\infty\le M_X$ almost surely. Then, we have
$$E\left[\left.\left|\frac1n\sum_{i=1}^n \epsilon_iX_i\right|_\infty\right|\{X_i\}_{i=1}^n\right]\le \sigma M_X\sqrt{\frac{2\log(2d)}{n}}.$$
\end{Lemma}
\begin{Proof} Take $v\in \R$ and $k\in\{1,\dots,d\}$. The random variable $n^{-1}\sum_{i=1}^n \epsilon_iX_{ki}$ follows a $\mathcal{N}(0,\sigma^2n^{-2}\sum_{i=1}^nX_{ki}^2)$ distribution conditional on $\{X_i\}_{i=1}^n$. Hence, by a standard property of Gaussian variables, $$E\left[\left.\exp((v/n)\sum_{i=1}^n \epsilon_iX_{ki})\right|\{X_i\}_{i=1}^n\right]\le \exp(v^2\sigma^2\sum_{i=1}^nX_{ki}^2/(2n^2))\le  \exp(v^2\sigma^2M_X^2/(2n)).$$
 For all $v>0$, this implies
\begin{align*}
&E\left[\left.\left|\frac1n \sum_{i=1}^n\epsilon_iX_{i}\right|_\infty\right|\{X_i\}_{i=1}^n\right]\\
&= vE\left[\left.\left|\frac1n \sum_{i=1}^n\epsilon_iX_{i}\right|_\infty/v\right|\{X_i\}_{i=1}^n\right]\\
&= v E\left[\left.\log\left(\exp\left(\left|\frac1n \sum_{i=1}^n\epsilon_iX_{i}\right|_\infty/v\right)\right)\right|\{X_i\}_{i=1}^n\right]\\
&\le v \log\left(E\left[\left.\exp\left(\left|\frac1n \sum_{i=1}^n\epsilon_iX_{i}\right|_\infty/v\right)\right|\{X_i\}_{i=1}^n\right]\right)\\
&\le v\log\left(\sum_{k=1}^d E\left[\left.\exp\left(\frac{1}{n} \sum_{i=1}^n\epsilon_iX_{ki}/v\right)\right|\{X_i\}_{i=1}^n\right] + E\left[\left.\exp\left(-\frac1n \sum_{i=1}^n\epsilon_iX_{ki}/v\right)\right|\{X_i\}_{i=1}^n\right]\right)\\
&\le v\log\left(2d\exp\left(\frac{\sigma^2M_X^2}{2nv^2}\right)\right)= v\left(\log(2d)+ \frac{\sigma^2M_X^2}{2nv^2}\right),
\end{align*}
where the first inequality is due to Jensen's inequality. Choosing $v=\sqrt{\sigma^2M_X^2/(2n\log(2d))}$, we obtain the result.
\end{Proof}

\section{Proof of Proposition 4}\label{sec.prop4}
To prove Proposition 4, it suffices to show that the assumptions of Theorem 2 hold and apply the latter theorem. To prove that Assumptions 2 and 4 are satisfied, we use Proposition 3, of which the conditions are shown to hold in the next subsection. Then, we prove that Assumptions 1 and 3 are satisfied in Section \ref{Lemmarobust}, thanks to Proposition 2. Finally, remark that $\lambda_n|\theta_0|_1\to 0$ since $\log(nd)|\theta_0|_1^2n^{-1}\to 0$ and $\lambda_n=2r_n$, where $r_n$ is defined in Proposition 3.

\subsection{Showing that the assumptions of Proposition 3 hold}
First, we show that Assumption 5 holds. Remark that Assumptions 5 (i) and (ii) are implied by Assumptions 6 and 7 (i), respectively. The Lipschitz property of $\ell$ is a direct consequence of $|\rho'(t)|\le M_\rho$ which is included in Assumption 7 (i).

It remains to prove that $\log(d)|B|_1^2n^{-1}\to 0$. Since $R_n(\theta_{0})\le M_\rho$, we have $|B|_1\le \lambda_n^{-1}M_\rho+|\theta_{0}|_1$. Hence, 
$$\log(d)|B|_1^2n^{-1}\le 2
\frac{\log(d)}{(32   M_X)^2 M_\rho  \log(4nd ) }    + 2 \log(d)|\theta_{0}|_1^2n^{-1} , $$
which goes to $0$ by assumption.

\subsection{Showing that the assumptions of Proposition 2 hold}\label{Lemmarobust}

To obtain both Assumptions 1 and 2 we only need to show that Condition (i) and (ii) in Proposition 2 are satisfied under Assumptions 6 and 7.

\paragraph{Proof of (i)  in Proposition 2.}

Take $\gamma>0$ and $\theta\in \Theta,\ \theta\ne \theta_{0},\ |\theta-\theta_{0}|_2^2\le \gamma$. For $s\ge 0$, define the event $$\mathcal{E}_s=\left\{|X^\top(\theta-\theta_{0})|\le s\right\}.$$
For $s>0$, let $L(s)=\inf_{0<|t|\le s}\frac{g(t)}{t}$.  Let us show that $L(s)>0$ for all $s>0$. Since, $\rho'$ is odd and $\epsilon$ has a symmetric distribution, $g$ is also odd and as a result $L(s) =\inf_{0<t\le s}\frac{g(t)}{t}$. Recall that we assumed that $g'(0)>0$. Hence, since
 $$g'(0)=\lim_{t\to0}\frac{g(t)-g(0)}{t}=\lim_{t\to0}\frac{g(t)}{t},$$ there exists $\eta>0$ such that $\inf_{0<t\le \eta}\frac{g(t)}{t}>0$. Next, since $t\mapsto g(t)/t$ is continuous and strictly positive on $[\eta,\infty)$, we have $\inf_{\eta<t\le s}\frac{g(t)}{t}>0$ for all $s>0$. This yields that $L(s)>0$ for all $s>0$.
As a result, for all $s>0$, we have
\begin{align*}
\nabla R(\theta)^\top(\theta-\theta_{0})&=E\left[E\left[ \rho'( X^\top(\theta-\theta_{0}) +\epsilon) | X  \right] X^\top(\theta-\theta_{0}) \right]\\
&=E[g(X^\top(\theta-\theta_{0}))X^\top(\theta-\theta_{0})]\\
&\ge E[g(X^\top(\theta-\theta_{0}))X^\top(\theta-\theta_{0})1_{\mathcal{E}_s}]\\
&\ge L(s) E[\{X^\top(\theta-\theta_{0})\}^21_{\mathcal{E}_s}],
\end{align*}
where, in the first inequality, we used the fact that $g(X^\top(\theta-\theta_{0}))X^\top(\theta-\theta_{0})$ is positive almost surely since $g$ is odd.
This shows that $\nabla R(\theta)^\top(\theta-\theta_{0})\ge 0$.

\paragraph{Proof of (ii)  in Proposition 2.}
Write
$$E[\{X^\top(\theta-\theta_{0})\}^21_{\mathcal{E}_s}]=E[\{X^\top(\theta-\theta_{0})\}^2] - E[\{X^\top(\theta-\theta_{0})\}^21_{\{X^\top(\theta-\theta_{0})\}^2 >s }].$$
Noting that $E[\{X^\top(\theta-\theta_{0})\}^2]  \geq \rho_X |\theta-\theta_{0}|_2^2 $, we now provide an upper  bound on $E[\{X^\top(\theta-\theta_{0})\}^21_{\{X^\top(\theta-\theta_{0})\}^2 >s }]$ when $s = s_\gamma$ with
\begin{align*}
&s_\gamma = 2 \sqrt { M_X^2 \gamma^2  \log(   8\sqrt 2  M_X^2 / \rho_X    ) }.
\end{align*}
Let  $ s_\theta = 2 \sqrt {  M_X^2|\theta-\theta_{0}|_2^2  \log(   8\sqrt 2  M_X^2 /  \rho_X ) }$. Since $X^\top(\theta-\theta_{0})$ is $v = M_X^2|\theta-\theta_{0}|_2^2$-sub-Gaussian, in virtue of Lemma \ref{lem:subG2} (with $\delta = \rho_X/(2M_X^2)$)  and the fact that $s_\theta\leq  s_\gamma$, we obtain 
\begin{align*}
 E[\{X^\top(\theta-\theta_{0})\}^21_{\{X^\top(\theta-\theta_{0})\}^2 >s_\gamma }]&\leq  E[\{X^\top(\theta-\theta_{0})\}^21_{\{X^\top(\theta-\theta_{0})\}^2 >s_\theta  }]   \\
 &\leq  |\theta-\theta_{0}|_2^2 \rho_X / 2
\end{align*}
Back to the initial decomposition, it follows that
$$E[\{X^\top(\theta-\theta_{0})\}^21_{\mathcal{E}_{s_\gamma}}] \geq  |\theta-\theta_{0}|_2^2 \rho_X /2.$$
and finally we showed that
$$\nabla R(\theta)^\top(\theta-\theta_{0})\ge L( s_{\gamma}) \frac{\rho_X}{2} |\theta-\theta_{0}|_2^2.$$
Hence, (ii) holds true with $c(\gamma)=L( s_{\gamma}) \frac{\rho_X}{2}$. The mapping $c(\cdot)$ is decreasing because $L(\cdot)$ is decreasing (by definition) and $ s_\gamma$ is increasing in $\gamma$.

\section{Proof of Proposition 5}\label{sec.prop5}
The proof of Proposition 5 leverages Propositions 2 and 3 in the same manner as the proof of Proposition 4.

\subsection{Showing that the assumptions of Proposition 3 hold}

Let us show that Assumption 5 holds. Assumptions 5 (i) and (ii) are a direct consequence of Assumptions 6 and 8 (i) and the fact that $Y$ is binary. 

Now, we show the Lipschitz property, corresponding to Assumption 5 (iii). For $t,t'\in \R$ and $y\in\{0,1\}$, we have
\begin{align*}
|\ell(t,y)-\ell(t',y)|&=|(y-\sigma(t))^2- (y-\sigma(t'))^2|\\
&=|2(y-\sigma(t'))(\sigma(t')-\sigma(t))+ (\sigma(t')-\sigma(t))^2|\\
&\le 3| \sigma(t')-\sigma(t) |\\
&\le 3M_\sigma  |t-t'| ,
\end{align*}
where we used the fact that the range of $\sigma$ belongs to $[0,1]$. This shows Assumption 5 (iii).

It  remains to prove that $\log(d)|B|_1^2n^{-1}\to 0$. Since $R(\theta_{0})\le 1$, we have $|B|_1\le \lambda_n^{-1}+|\theta_{0}|_1$. By  Assumption 8 (iii), we have $|\theta_{0}|_1\le \sqrt{s_0} |\theta_{0}|_2\le  \sqrt{s_0} M_0$. Hence, 
$$\log(d)|B|_1^2n^{-1}\le 2
\frac{\log(d)}{(96   M_X M_\sigma)^2   \log(4nd ) }    + 2 \log(d)|\theta_{0}|_1^2n^{-1} , $$
which goes to $0$ by assumption.


\subsection{Showing that the assumptions of Proposition 2 hold}\label{E2march}

\paragraph{Proof of (i)  in Proposition 2.} Taking the derivative inside the expectation, we obtain
\begin{align*}
 \nabla R(\theta )  =  2 \mathbb E [ - \sigma' ( X^\top\theta) X  (  Y - \sigma ( X^\top\theta)  )  ],
\end{align*}
and using that $ \mathbb E [Y| X] = \mathbb P ( Y=1 | X)  =  \sigma(X^\top\theta_{0} )  $, it follows that
\begin{align*}
 \nabla R(\theta )  =  2 \mathbb E [  \sigma' ( X^\top\theta) X  (    \sigma ( X^\top\theta)  -  \sigma(X^\top\theta_{0} )    )  ].
\end{align*}
Now introduce $\theta _ t  = t \theta + (1-t) \theta_{0}$ and $F(t) = \sigma ( \theta_t^\top X)$, $0\leq t\leq 1$, we can write 
\begin{align*}
\sigma ( X^\top\theta) -   \sigma(X^\top\theta_{0} ) &= F(1) - F(0) \\
&= \int _0^1 F'(t) dt\\ 
& = \{(\theta - \theta_{0}) ^\top X\} \int_0^1 \sigma' (\theta_t ^\top X) dt 
\end{align*}
in order to obtain
\begin{align*}
 \nabla R(\theta )^\top (\theta - \theta_{0})  & = 2\mathbb E [ L (X,\theta)  \{ X^\top  (\theta - \theta_{0})\} ^2 ]
\end{align*}
with  $ L (X,\theta) =  \sigma' ( X^\top\theta) \int_0 ^1 \sigma' (\theta _t ^\top X ) dt$. It implies that (i) holds by Assumption 7 (ii).

\paragraph{Proof of (ii)  in Proposition 2.} 
Introduce the two events
\begin{align*}
&\mathcal E _s = |X^\top  (\theta - \theta_{0}) |\leq s \qquad \text{and} \qquad  \mathcal F _u = |X^\top   \theta_{0} |\leq u.
\end{align*}
If both are realized, then $ | \theta_tX| \leq  s+ u $ for all $t\in [0,1]$. Consequently, 
$$L (X,\theta) 1_{\mathcal E _s}  1_{\mathcal F _{u}  }\geq L(s,u) : =  \inf_{|t| \leq s+u }    \sigma' ( t )  ^2 ,$$
and therefore we find
\begin{align*}
 \nabla R(\theta )^\top (\theta - \theta_{0}) & \geq  2 L(s,u ) \mathbb E [ 1_{\mathcal E _s}  1_{\mathcal F_{u}  }  \{X^\top  (\theta - \theta_{0})\}^2 ].
 \end{align*}
Recalling that $ \rho_X$ is a lower bound on the smallest eigenvalue of $ \mathbb E [XX^\top ]$,  it follows that
 \begin{align*}
  \nabla R(\theta )^\top (\theta - \theta_{0})  & \geq  2 L(s,u ) ( \rho_X | \theta - \theta_{0} |_2^2  - R(s,s_0))
\end{align*}
with $ R(s,u ) = \mathbb E [ 1_{\mathcal E _s^c \cup \mathcal F_u^c  }  \{X^\top  (\theta - \theta_{0})\}^2 ]$. The union bound gives
\begin{align*}
 R(s,u ) \leq  \mathbb E [ (1_{\mathcal E _s^c } + 1_{  \mathcal F_u^c  }   )\{X^\top  (\theta - \theta_{0})\}^2 ].
\end{align*}
Let
\begin{align*}
&s_\theta =  2 \sqrt { | \theta - \theta_{0}|^2_2 M_X^2   \log(   16\sqrt 2 M_X^2 / \rho_X ) }  ,  & u_0 =  2 \sqrt { |  \theta_{0}|^2_2 M_X^2   \log(   16 \sqrt 2 M_X^2  / \rho_X ) }  .
\end{align*}
In virtue of Lemma \ref{lem:subG2}, we have
\begin{align*}
R(s_\theta , u_0 ) \leq \frac {\rho_X}  {4} |\theta-\theta_{0}|_2^2  +  \frac {\rho_X}  {4}  |\theta-\theta_{0}|_2^2   =  \frac {\rho_X}  {2}  |\theta-\theta_{0}|_2^2 
\end{align*}
implying that for all $\theta \in \mathbb R^d$,
\begin{align*}
 \nabla R(\theta )^\top (\theta - \theta_{0})   &\geq   L(s_\theta ,u_0 )  \rho_X | \theta - \theta_{0} |_2^2.
\end{align*}
Now take the infimum and use that $ s\mapsto  L(s ,u_0 ) $ and $ u \mapsto  L(s ,u ) $ are  nonincreasing functions to get
\begin{align*}
\inf_{\theta\in \Theta, \, |\theta- \theta_{0}|_2\leq \gamma } \frac{  \nabla R(\theta )^\top (\theta - \theta_{0})  }{| \theta - \theta_{0} |_2^2} & \geq   L(s_\gamma ,u_0 )  \rho_X  \geq   L(s_{\gamma,0} ,s_{\gamma,0}  )  \rho_X
\end{align*} 
with $ s_\gamma  = 2 \sqrt { \gamma ^2  M_X^2   \log(   16\sqrt 2 M_X^2 / \rho_X ) }  $ and $s_{\gamma,0}  =2 M_X (\gamma \vee M_0)  \sqrt {  \log(   16\sqrt 2 M_X^2 / \rho_X ) } $.

\section{Proof of Proposition 6}\label{sec.prop6}

To prove Proposition 6, we show that the assumptions of Theorem 2 hold and apply the latter theorem. 
\subsection{Proof that Assumptions 2 and 4 hold} 
Using $\epsilon|X\sim\mathcal{N}(0,\sigma^2)$, we obtain
 $$\nabla R(\theta )  =  2 \mathbb E [  f' ( X^\top\theta) X  (    f ( X^\top\theta)  -  f(X^\top\theta_{0} )    )  ].$$ Hence, $\nabla R(\theta)$ in the application of nonlinear least squares has the same shape as in the application of binary classification (albeit with $\sigma$ replaced by $f$). Since we make similar assumptions on $\sigma$ in Section 4.2 and on $f$ in Section 4.3, Assumptions 1 and 3 can be proved to hold exactly in the same manner as in the proof of Proposition 5 in Section \ref{E2march}. The proof is therefore omitted.
\subsection{Proof that Assumptions 2 and 4 are satisfied} 
Using the decomposition of equation (7) in the main text, we have, for all $\theta\in  \Theta$, \begin{align*}R(\theta)&= R_{1}(\theta)-2R_{2}(\theta)+E[\epsilon^2];\\
\widehat{R}(\theta)&= \widehat{R}_{1}(\theta)-2\widehat{R}_{2}(\theta)+\frac1n\sum_{i=1}^n\epsilon_i^2,
\end{align*} where 
\begin{align*}
&R_{1}(\theta) = E[(f(X^\top\theta)-f(X^\top\theta_{0}))^2]; \\
& \widehat R_{1}(\theta)=\frac1n\sum_{i=1}^n (f(X_i^\top\theta)-f(X_i^\top\theta_{0}))^2;\\
 & \widehat R_{2}(\theta)= \frac1n\sum_{i=1}^n\epsilon_i(f(X_i^\top\theta)-f(X_i^\top\theta_{0})).
\end{align*}
For $\theta\in \Theta,$ let $\widehat{\Delta}(\theta)=\widehat R(\theta)- R(\theta)$, $\widehat{\Delta}_{1}(\theta)=\widehat{R}_{1}(\theta)-R_{1}(\theta)$ and  $\widehat{\Delta}_{2}(\theta)=\widehat{R}_{2}(\theta)$ (remark that $  E [ \widehat R_2 (\theta ) ] = 0$). Lemmas \ref{R1} and \ref{R22} in the next two subsections show that Assumptions 2 and 4 hold for $\widehat  R_{1}$ and $\widehat R_{2}$, respectively. 

As a result, we have 
\begin{align*}\sup_{\theta\in B} |\widehat \Delta(\theta)|
&\le \sup_{\theta\in B} |\widehat \Delta_{1}(\theta)|+2\sup_{\theta\in B} |\widehat \Delta_{2}(\theta)|+ \left|\frac{1}{n}\sum_{i=1}^n \epsilon_i^2-E[\epsilon^2]\right|=o_P(1),
\end{align*}
where we used the triangle inequality, Lemmas \ref{R1} (i) and \ref{R22} (i) and the fact that $\left|n^{-1}\sum_{i=1}^n \epsilon_i^2-E[\epsilon^2]\right|=o_P(1)$ by the law of large numbers. This shows that Assumption 2 is satisfied.

Next, let $\delta_n =\sqrt{\log(2d)/n}$ and $r_n=5(K_2\vee K_3) \sqrt{\log(4nd)/n}$, where $K_2$ and $K_3$ are defined in Lemmas \ref{R1} (ii) and \ref{R22} (ii). Remark that, for all $\theta \in B$, we have \begin{equation}\label{decompnls2}\widehat \Delta(\theta) - \widehat \Delta(\theta_{0})=\widehat \Delta_{1}(\theta) - \widehat \Delta_{1}(\theta_{0}) - 2 (\widehat \Delta_{2}(\theta) - \widehat \Delta_{2}(\theta_{0})) .\end{equation} It holds that

\begin{align*}&\P\left(\sup_{\theta\in B} \frac{|\widehat \Delta(\theta)- \widehat \Delta(\theta_{0})|}{|\theta-\theta_{0}|_1\vee \delta_n} >  r_n\right)\\
&\le \P\left(\sup_{\theta\in B} \frac{|\widehat \Delta_{1}(\theta)- \widehat \Delta_{1}(\theta_{0})|}{|\theta-\theta_{0}|_1\vee \delta_n} +2\sup_{\theta\in B} \frac{|\widehat \Delta_{2}(\theta)- \widehat \Delta_{2}(\theta_{0})|}{|\theta-\theta_{0}|_1\vee \delta_n} > r_n\right)\\
&\le \P\left(\sup_{\theta\in B} \frac{|\widehat \Delta_{1}(\theta)- \widehat \Delta_{1}(\theta_{0})|}{|\theta-\theta_{0}|_1\vee \delta_n} >  \frac{r_n}{2}\right)+ \P\left(\sup_{\theta\in B} \frac{|\widehat \Delta_{2}(\theta)- \widehat \Delta_{2}(\theta_{0})|}{|\theta-\theta_{0}|_1\vee \delta_n} > \frac{r_n}{4}\right)\to 0,
\end{align*}
where in the first inequality, we used \eqref{decompnls2} and the triangle inequality, and, in the second inequality, we leveraged the union bound. The fact that the limit of the probability is $0$ is a consequence of Lemmas \ref{R1} (ii) and \ref{R22} (ii). This proves that Assumption 4 is satisfied and concludes the proof.
\subsection{On the term $\widehat{R}_{1}$}\label{usefullmm}
Recall that, for $\theta\in \Theta,$  $\widehat{\Delta}_{1}(\theta)=\widehat{R}_{1}(\theta)-R_{1}(\theta)$.
\begin{Lemma}\label{R1}
 Under Assumptions 6 and 9, we have 
\begin{itemize}
\item[(i)] $\sup_{\theta\in B}\left| \widehat{\Delta}_{1}(\theta)\right|=o_P\left(1\right).$
\item[(ii)] For all $\eta>0$,
\begin{align*}
 \P\left(\sup_{\theta\in B}\frac{\left|\widehat{\Delta}_{1}( \theta)  -  \widehat{\Delta }_{1}( \theta_{0})   \right|}{|\theta-\theta_{0}|_1\vee \delta_n}   \leq   r_n\right)\to 1,
\end{align*}
 where 
$$\delta_n =\sqrt{\frac{\log(2d)}{n}},\ r_n=K_1\sqrt{\frac{\log(4nd)}{n}},$$ for some constant $K_1$.
\end{itemize}
\end{Lemma}
\begin{Proof}
Let $\mathcal{T}$ be defined as in Section 3.2 and $\widetilde{\mathcal{Y}}= \{x^\top\theta_{0}:\ x\in\mathcal{X}_n\}$. For $\theta\in \Theta$, we have $$ R_{1}(\theta)=E[\ell(X^\top\theta, \widetilde Y)],$$
where $\ell:(s,y)\in \mathcal{T}\times \widetilde{\mathcal{Y}}\mapsto (f(s)-f(y))^2$ and $\widetilde Y=X^\top\theta_{0}$. We also let 
$$\widehat{ R}_1(\theta)=\frac1n\sum_{i=1}^n\ell(X_i^\top\theta, \widetilde Y_i),$$
with $\widetilde Y_i=X_i^\top\theta_0$. 
With these notations, $R_{1}$ belongs to the class of true risks $R$ considered in Section 3.2. It is clear that Assumption 5 (i) and (ii) holds thanks to Assumptions 6 and 9 (i), respectively. Moreover, Assumption 5 (iii) is satisfied because for $t,t'\in \mathcal{T}$ and $y\in\widetilde{\mathcal{Y}}$, we have
\begin{align*}
|\ell(t,y)-\ell(t',y)|&=|(f(t)-f(y))^2- (f(t')-f(y))^2|\\
&= | (f(t)- f(t')) ( f(t) + f(t') - 2 f(y) )|\\
&\le 4M_f|f(t')-f(t)|\\
&\le 4M_f^2 |t-t'|.
\end{align*}
We can show that $\log(d)|B|_1^2n^{-1}\to 0$ as in the proof of Proposition 5.
The result of the Lemma is then a direct consequence of Proposition 3 applied to $ R_{1}$ and $\widehat R_{1}$.
\end{Proof}

\subsection{On the term $\widehat{R}_{2}$}
Recall that, for $\theta\in \Theta$,  $\widehat{\Delta}_{2}(\theta)=\widehat{R}_{2}(\theta)$.
\begin{Lemma}\label{R21}
For a nonempty bounded set $C \subset \Theta$, define $\mu_{C}= \sup_{\theta\in C} |\theta-\theta_{0}|_1$ and 
$$D =\sup\limits_{\theta\in C}\left|\widehat{\Delta}_{2}(\theta)-\widehat{\Delta}_{2}(\theta_{0})\right|.$$
 Under Assumptions 6 and 9, we have 
$$\P\left(D \ge  2K_2\sqrt{\frac{2\log(2d)}{n}}\mu_C+K_2\mu_Ct\right)\le 2\exp\left(-\frac{nt^2}{2}\right),$$for some constant $K_2>0$.
\end{Lemma}
\begin{Proof} 
Without loss of generality, we can assume that $\theta_0\in C$. If not it suffices to replace $C$ by $C\cup \{\theta_0\}$ and check that nothing changes in the statement of the lemma. 
Notice that  $ \widehat{R}_{2}(\theta_{0})=0$. Then, $\widehat{\Delta}_{2}(\theta)-\widehat{\Delta}_{2}(\theta_{0})=\widehat{R}_{2}(\theta)$ and $D=\sup_{\theta\in C}\left|\widehat{R}_{2}(\theta)\right|.$ The proof follows from the fact that $\{\widehat R_{2}(\theta)\}_{\theta\in C}$ is a Gaussian process given $\{X_i\}_{i=1}^n $. \\

\noindent\textit{Bounding $E[D]$.} 
Since $\epsilon| X\sim\mathcal{N}(0,\sigma^2)$, for $\theta,\gamma\in \Theta$, we have
$$\widehat R_{2}(\theta)-\widehat R_{2}(\gamma)|\{X_i\}_{i=1}^n\sim \mathcal{N}\left(0,\rho^2(\theta,\gamma)\right),$$
where $\rho(\theta, \gamma)=n^{-1} \sigma  \sqrt{\sum_{i=1}^n(f(X_i^\top\theta)-f(X_i^\top\gamma))^2}.$
Hence, by the Gaussian concentration inequality (Theorem 5.6 in \cite{boucheron2013concentration}), we have, for all $t\ge 0$, 
$$\P\left(\left.|\widehat R_{2}(\theta)-\widehat R_{2}(\gamma)|\ge t\rho(\theta, \gamma)\right|\{X_i\}_{i=1}^n\right)\le e^{-\frac{t^2}{2}}.$$
Remark that, by Assumption 9 (ii), $\rho^2(\theta, \gamma)\le  M_f^2 \sigma  ^2 n^{-2}\sum_{i=1}^n(X_i^\top(\theta-\gamma))^2.$
This yields, 
$$\P\left(\left.|\widehat R_{2}(\theta)-\widehat R_{2}(\gamma)|\ge t  M_f^2 \sigma  ^2 n^{-2}\sum_{i=1}^n(X_i^\top(\theta-\gamma))^2\right|\{X_i\}_{i=1}^n\right)\le e^{-\frac{t^2}{2}}.$$
As a result, conditional on $\{X_i\}_{i=1}^n$, the process $\{\widehat R_{2}(\theta)\}_{\theta\in C}$ has the same tail bound as the Gaussian process
$\{P_\theta=M_f \sigma n^{-1}\sum_{i=1}^n \epsilon_iX_i^\top(\theta-\theta_{0})\}_{\theta\in C}.$ Hence, using Theorem 2.1.5 in \cite{talagrand2005generic} and the fact that $\theta_0\in C$, we get that there exists a constant $L>0$ such that 
$$E[D|\{X_i\}_{i=1}^n]\le E\left[\left.\sup_{\theta,\gamma\in C}|\widehat R_{2}(\theta)-\widehat R_{2}(\gamma)|\right|\{X_i\}_{i=1}^n\right]\le LE\left[\left.\sup_{\theta\in C} |P_\theta|\right|\{X_i\}_{i=1}^n\right].$$
Next, we have 
$$E\left[\left.\sup_{\theta\in C} |P_\theta|\right|\{X_i\}_{i=1}^n\right]\le M_f\sigma   E\left[\left.\left|\frac1n\sum_{i=1}^n\epsilon_iX_i\right|_\infty\right|\{X_i\}_{i=1}^n\right] \mu_C \le \sigma M_fM_X\sqrt{\frac{2\log(2d)}{n}}\mu_C,$$
where the last inequality is due to Lemma \ref{tech2}. By the law of iterated expectations, we get 
\begin{equation}\label{boundnls}
E[D]\le \sigma L M_fM_X\sqrt{\frac{2\log(2d)}{n}}\mu_C.
\end{equation}
\\

\noindent\textit{Concentrating $D$ around $E[D]$.} 
We have $$\widehat R_{2}(\theta)|\{X_i\}_{i=1}^n\sim \mathcal{N}\left(0,\rho^2(\theta,\theta_{0})\right),$$ By the Gaussian bound for processes corresponding to Theorem 5.8 in \cite{boucheron2013concentration}, it holds that, for all $t\ge 0$, 
\begin{equation}\label{Gaussian bound} \P\left(|D-E[D]|\ge t\bar\rho)|\{X_i\}_{i=1}^n\right)\le 2e^{-\frac{t^2}{2}},\end{equation}
where $\bar \rho= {\sup_{\theta\in C} |\rho(\theta,\theta_{0}))|}.$
Remark that, for all $\theta\in C$, it holds that
$$\rho^2(\theta,\theta_{0})\le n^{-2} \sigma ^2 M_f^2\sum_{i=1}^n(X_i^\top(\theta-\theta_{0}))^2\le n^{-1}  \sigma  ^ 2 M_f^2M_X^2|\theta-\theta_{0}|_1\le n^{-1}  \sigma  ^ 2  M_f^2M_X^2\mu_C,$$
almost surely, where the second inequality leverages Assumption 6 and Hölder's inequality.
As a result, we have $\bar \rho^2\le   \sigma  ^ 2  M_f^2M_X^2\mu_C^2n^{-1}$ almost surely. By \eqref{Gaussian bound} and the bound on $E[D]$ in \eqref{boundnls}, this yields the result the bound on the probability given $ \{X_i\}_{i=1}^n$. The unconditional bound can be directly obtained through the law of iterated expectations.
\end{Proof}

\begin{Lemma}\label{R22}
 Under Assumptions 6 and 9, we have 
\begin{itemize}
\item[(i)] $\sup_{\theta\in B}\left| \widehat{\Delta}_{2}(\theta)\right|=o_P\left(1\right).$
\item[(ii)] For all $\eta>0$,
\begin{align*}
 \P\left(\sup_{\theta\in B}\frac{\left| \widehat{\Delta}_{2}(\theta)-  \widehat{\Delta}_{2}(\theta_{0}) \right|}{|\theta-\theta_{0}|_1\vee \delta_n}   \leq   r_n\right)\to 1,
\end{align*}where 
$$\delta_n =\sqrt{\frac{\log(2d)}{n}},\ r_n=K_3\sqrt{\frac{\log(4nd)}{n}},$$ where $K_3>0$ is a constant.
\end{itemize}
\end{Lemma}
\begin{Proof} Lemma \ref{R21} is a result similar to Lemma 2 in the main text. So, the proof can proceed similarly to the proof of Proposition 3 in Sections C.2.2 and C.2.3 of the main text, and is, therefore, omitted.
\end{Proof}

\bibliographystyle{apalike}
\bibliography{Supplement}